\documentstyle{amsppt}
\magnification=\magstep1
\pageheight{185mm}
\pagewidth{138mm}
\topskip= 1.25cm
\nologo
\parindent=1em
\parskip=4pt
\baselineskip=12pt
\loadbold
\loadeufm

\topmatter
\title
Pro-$p$ link groups and $p$-homology groups
\endtitle
\author
Jonathan Hillman$^{1}$, Daniel Matei$^{2}$ and Masanori Morishita$^{3}$
\endauthor
\leftheadtext{Hillman, Matei and Morishita}
\rightheadtext{Pro-$p$ link groups and $p$-homology groups}
\endtopmatter
{\sevenrm
\roster
\leftline{1.  School of Mathematics and Statistics, The University of
Sydney,  Sydney, NSW 2006, Australia}
\leftline{2. Graduate School of Mathematical Sciences, University of
Tokyo,  3-8-1 Komaba, Meguro, Tokyo 153-8914, Japan}
\leftline{3.  Department of Mathematics, Kanazawa University,
Kakuma-machi, Kanazawa, Ishikawa, 920-1192, Japan}
\endroster}

{\eightrm
\leftline{\bf Contents}
\leftline{Introduction}
\leftline{1. Pro-$p$ completion of a link group}
\leftline{2. $p$-adic Milnor invariants }
\leftline{3. Completed Alexander modules}
\leftline{4. Galois module structure of the $p$-homology group of a
$p$-fold cyclic branched cover }
\leftline{5. Iwasawa type formulas for the $p$-homology groups of
$p^{m}$-fold cyclic branched covers}}
\medskip
\document

\head {\bf Introduction } \endhead

Let $L$ be a tame link in the 3-sphere $S^{3}$ consisting of $n$  knots
$K_{1}, \cdots, K_{n}$ and let $G_{L}$ be the link group
$\pi_{1}(X_{L}), X_{L} =  S^{3} \setminus L$. For a prime number $p$,
let  $\widehat{G_{L}}$ denote the pro-$p$  completion of the group
$G_{L}$, $\widehat{G_{L}} = \varprojlim G_{L}/N$ where $N$ runs over
normal subgroups of $G_{L}$ having $p$-power indices. By a theorem of J.
Milnor [Mi], it is shown  that $\widehat{G_{L}}$ has the following
simple presentation as a pro-$p$ group
$$ \widehat{G_{L}} = \langle x_{1},\cdots ,x_{n} \; | \;
[x_{1},y_{1}]=\cdots =[x_{n},y_{n}] = 1 \rangle$$
where $x_{i}$ and $y_{i}$ represent the meridian and longitude around
$K_{i}$ respectively (Theorem 1.2.1). The purpose of this paper is to
use the pro-$p$ link group $\widehat{G_{L}}$ and the associated
group-theoretic invariants for the study of the $p$-homology groups of
$p^{m}$-fold cyclic branched covers of $S^{3}$ along $L$, following the
analogies between link theory and  number theory [Mo1$\sim$4],[Rez1,2].
The invariants we derive from $\widehat{G_{L}}$ are the $p$-adic Milnor
invariants and the completed Alexander module over the formal power
series ring $\widehat{\Lambda_{n}} = {\bold Z}_{p}[[X_{1},\cdots
,X_{n}]]$ with coefficients in the ring ${\bold Z}_{p}$ of $p$-adic
integers. The tool involved here is the Fox differential calculus on a
free pro-$p$ group [Ih]. Although these invariants are simply $p$-adic
analogues of the usual Milnor invariants and Alexander modules,  it is
natural to work over  $\widehat{\Lambda_{n}}$ since the completed
Alexander module can be presented over $\widehat{\Lambda_{n}}$ by a sort
of {\it universal $p$-adic higher linking matrix} $\widehat{T_{L}}$,
called the {\it $p$-adic Traldi matrix}. This is  defined in terms of
the $p$-adic Milnor numbers and we can derive from $\widehat{T_{L}}$
 systematically the ``$p$-primary" information on the homology of $p^{m}$-fold 
 branched covers of $L$. This is an idea analogous to
Iwasawa theory [Iw] which may also be regared as a $p$-adic
strengthening of the method employed by W. Massey [Mas] and L. Traldi
[T]. We note that the method using the truncated Traldi matrices was
considered in [Mat] to study the homology of unbranched covers.
\medskip

The homology of cyclic branched covers of a link $L$ is one of the
basic invariants of $L$ and has been extensively investigated by many
authors. The Betti number and the order have been determined in terms of
the Alexander (Hosokawa) polynomial ([HK],[MM],[S1] etc) and further the
(Galois) module structure has been studied ([Da],[HS],[S2] etc), however
most results are concerned mainly with the part which is prime to the
covering degree. In [Rez1,2], A. Reznikov studied the $p$-homology of
$p$-fold branched covers after the model of the classical problem on
$p$-ideal class groups in number theory (see also [Mo1]). In this paper,
we push this line of study in {\it arithmetic topology} further and
determine the Galois module structure of the $p$-homology of a $p$-fold
branched cover along a link completely in terms of the $p$-adic higher
linking matrices. To be precise, let $M$ be the $p$-fold cyclic branched
cover of $S^{3}$ along $L$ obtained from the completion of the $p$-fold
total linking cover of $X_{L}$ and let $\sigma$ denote a generator of
the Galois group of $M$ over $S^{3}$. The homology group $H_{1}(M,{\bold
Z}_{p}) = H_{1}(M,{\bold Z})\otimes_{\bold Z}{\bold Z}_{p}$ is then a
module over the complete discrete valuation ring $\widehat{\Cal O} :=
{\bold Z}_{p}[\langle \sigma \rangle]/(\sigma^{p-1}+\cdots +\sigma + 1)
= {\bold Z}_{p}[\zeta]$, $\zeta := \sigma$ mod $(\sigma^{p-1}+\cdots
+\sigma + 1)$. Assume that $H_{1}(M,{\bold Z})$ is finite. Then the
$p$-primary part  $H_{1}(M,{\bold Z}_{p})$ has $\frak{p}$-rank $n-1$
 ([Mo1],[Rez2]) so that it has  form
$$ H_{1}(M,{\bold Z}_{p}) = \bigoplus_{i=1}^{n-1}\widehat{\Cal
O}/\frak{p}^{a_{i}} \;\; (a_{i} \geq 1)$$
as $\widehat{\Cal O}$-module where $\frak{p} :=(\zeta -1)$ is the
maximal ideal of $\widehat{\Cal O}$. Hence the determination of the
Galois module structure of $H_{1}(N,{\bold Z}_{p})$ is equivalent to
that of the $\frak{p}^{k}$-rank
$$ e_{k} := \# \{ i \; | \; a_{i} \geq k \} \;\; (k\geq 1). $$
Our main result is to give formulas for $e_{k}$'s in terms of the
higher linking matrices obtained by specializing the truncated $p$-adic
Traldi matrices at $X_{1}=\cdots = T_{n} = \zeta -1$ (Theorem 4.1.3).
For the simplest case of $k=2$, our formula reads
$$  e_{2} = n-1 - \text{rank}_{{\bold F}_{p}}(C \; \text{mod} \; p)$$
where $C =(C_{ij})$ is the linking matrix defined  by $C_{ii} =
-\sum_{j\neq i}\text{lk}(K_{i},K_{j})$ and $C_{ij} =
\text{lk}(K_{i},K_{j})$ for $i\neq j$.
In view of the analogy between the linking number and the power residue
symbols  [Mo2,3], this is seen as a link-theoretic analog of L.
R\'{e}dei's formula for the $4$-rank of the class group of a quadratic
field ([R\'{e}d1]), and our general result was partly suggested by the
relation between R\'{e}dei's triple symbol and the $8$-rank of a class
group [R\'{e}d2]. In fact, the whole argument here can be translated
into arithmetic [Mo5]. In the last section, we study the asymptotic
behavior of the order $|H_{1}(M_{m},{\bold Z}_{p})|$ for the
$p^{m}$-fold cyclic branched cover $M_{m}$ as $m \rightarrow \infty$,
following Iwasawa theory on ${\bold Z}_{p}$-extensions [Iw]. Though our
results obtained in this paper are rather elementary, they seem to
indicate further possibilities of our arithmetic approach to link
theory.
\medskip
\flushpar
{\it Acknowledgement.} We would like to thank K. Murasugi for helpful
communication.
\medskip
\flushpar
{\it Notation.} Throughout this paper, we fix a  prime number $p$. We
denote by ${\bold F}_{p}$ the field with $p$ elements and by ${\bold
Z}_{p}$ the ring of $p$-adic integers. Let $\text{ord}_{p}$ denote the
additive $p$-adic valuation extended on the algebraic closure
$\overline{{\bold Q}_{p}}$ of the $p$-adic field ${\bold Q}_{p}$ with
$\text{ord}_{p}(p) = 1$ and set $|x|_{p} = p^{-\text{ord}_{p}(x)}, \; x
\in \overline{{\bold Q}_{p}}$. We use the letter $q$ to denote $p$ or
$0$. For a topological (possibly discrete) group $G$, we denote by
$G^{(k,q)}$ the $k$-th term of lower central $q$-series defined by
$G^{(1,q)} = G, G^{(k+1,q)} = (G^{(k,q)})^{q}[G^{(k,q)},G]$ where for
closed subgroups $A, B$ of $G$, $[A,B]$ stands for the closed subgroup
of $G$ generated by $[a,b] = aba^{-1}b^{-1}, a \in A, b \in B$. We
simply write $G^{(k)}$ for $G^{(k,0)}$, the $k$-th term of the lower
central series of $G$. For a pro-finite group $G$ and a complete local
ring $R$, we denote by $R[[G]]$ the completed group ring of $G$ over $R$
 [Ko,$\S 7$].
\medskip

\head{\bf 1. Pro-$p$ completion of a link group} \endhead
\medskip
\flushpar
1.1. {\it The pro-p completion of a link group.} Let $L$ be a tame link
in the 3-sphere $S^{3}$ consisting of $n$ component knots $K_{1}, \cdots
, K_{n}$ and let $G_{L}$ be the link group $\pi_{1}(X_{L})$, $X_{L} =
S^{3} \setminus L$.  After the work of K.T. Chen, J. Milnor [Mi]
derived the following information about the presentation of the
nilpotent quotient $G_{L}/G_{L}^{(k, q)}$. Let $F$ be the free group on
the $n$ words $x_{1}, \cdots , x_{n}$ where $x_{i}$ represents the
meridian $m_{i}$ around $K_{i}$ and let $\pi : F \rightarrow G_{L}$ be
the meridianal homomorphism defined by $\pi (x_{i}) = m_{i}$ $(1\leq i
\leq n)$.
\medskip
\flushpar
{\bf Theorem 1.1.1} ([Mi]). {\it For each $k \geq 1$ and $i$ ($1\leq
i\leq n$), there is a word $y_{i}^{(k)}$ in $x_{1},\cdots ,x_{n}$
 representing the image of the $i$-th longitude in the quotient
$G_{L}/G_{L}^{(k,q)}$ such that
$$y_{i}^{(k)} \equiv y_{i}^{(k+1)}\;\;\text{mod} \;\; F^{(k,q)} \tag
1.1.2 $$
 and such that $\pi : F \rightarrow G_{L}$ induces the isomorphism
$$ F/N_{k}F^{(k,q)} \; \simeq \;  G_{L}/G_{L}^{(k,q)} \tag 1.1.3 $$
where $N_{k}$ is the subgroup of $F$ generated normally by
$[x_{1},y_{1}^{(k)}], \cdots , [x_{n},y_{n}^{(k)}]$.}
\medskip

Let $\widehat{G_{L}}$ be the pro-$p$ completion of $G_{L}$, namely the
inverse limit $\varprojlim G_{L}/N$ of the tower of quotients $G_{L}/N$
which are finite $p$-groups. Since the quotients by the lower central
$p$-series of $G_{L}$ are cofinal in this tower, we have
$$ \widehat{G_{L}} \; = \; \varprojlim\limits_{k} G_{L}/G_{L}^{(k,p)}.
$$
Since $\{y_{i}^{(k)}F^{(k,p)}\}_{k \geq 1}$ forms an inverse system in
the tower $\{ F/F^{(k,p)}\}_{k \geq 1}$ by (1.1.2), we define the
pro-$p$ word $y_{i}$ to be $(y_{i}^{(k)}F^{(k,p)})$ in the free pro-$p$
 group $\widehat{F} := \varprojlim F/F^{(k,p)}$ which represents the
$i$-th ``longitude" in $\widehat{G_{L}}$ under the map $\widehat{\pi} :
\widehat{F}\rightarrow \widehat{G_{L}}$ induced by $\pi$. By taking the
inverse limit with respect to $k$ in the isomorphism (1.1.3) of finite
$p$-groups for $q = p$, we have the following
\medskip
\flushpar
{\bf Theorem 1.1.2.} {\it The map $\widehat{\pi}$ induces the
isomorphism of pro-$p$ groups
$$ \widehat{F}/\widehat{N} \; \simeq \; \widehat{G_{L}} $$
where $\widehat{N}$ is the closed subgroup of $\widehat{F}$ generated
normally by $[x_{1},y_{1}], \cdots ,[x_{n},y_{n}]$. In particular, we
have $\widehat{G_{K}} \; \simeq \; {\bold Z}_{p}$ for a knot $K$.}
\medskip
\flushpar
{\bf Remark 1.1.3.} (1) By the construction above, we note  $y_{i}
\equiv y_{i}^{(k)}$ mod $\widehat{F}^{(k,q)}$ ($F$ is embedded in
$\widehat{F}$).
\flushpar
(2) In view of the analogy between knots and primes, the pro-$p$ link
group $\widehat{G_{L}}$ is regarded as an analog of the maximal pro-$p$
 Galois group over  the rational number field ${\bold Q}$ unramified
outside prime numbers $p_{1},\cdots ,p_{n}$, $p_{i} \equiv 1$ mod $p$
 [Mo2].
\medskip

 Theorem 1.1.2 tells us that from the group-theoretic point of view, any
link of $n$ components looks like a pure braid link with $n$ strings
after the pro-$p$ completion. In particular, by applying the method of
D. Anick [A] to determine the graded quotients of the lower central
series of a pure braid link group to our pro-$p$ link group
$\widehat{G_{L}}$, we see that the pro-$p$ analog of Murasugi's
conjecture holds (cf. [L]).  We define the {\it mod $p$ linking diagram}
of $L$  to be the graph with vertices the components of $L$ and an edge
joining $K_{i}$ and $K_{j}$ if and only if the linking number
$\text{lk}(K_{i},K_{j}) \not\equiv 0$ mod $p$.
\medskip
\flushpar
{\bf Theorem 1.1.4.} {\it If the mod $p$ linking diagram of $L$ is
connected,
we have the isomorphisms}
$$  \widehat{G_{L}}^{(q)}/\widehat{G_{L}}^{(q+1)} \; \simeq  \;
\widehat{F}_{1}^{(q)}/\widehat{F}_{1}^{(q+1)} \times
\widehat{F}_{n-1}^{(q)}/\widehat{F}_{n-1}^{(q+1)}\;\; for \;\; q \geq
1,$$
{\it where $\widehat{F}_{r}$ denotes the free pro-$p$ group of rank
$r$.}
\medskip
\flushpar
1.2. {\it The $p$-goodness of a link group.} Let $G$ be a group and
$\widehat{G}$ be the pro-$p$ completion of $G$. We then call $G$ {\it
$p$-good} if the natural map $G \rightarrow \widehat{G}$ induces the
isomorphisms on cohomology $ H^{q}(\widehat{G},M) \; \overset \sim \to
\rightarrow \; H^{q}(G,M)$  for all $q \geq 0$ for any finite
$p$-primary $\widehat{G}$-modules $M$ (cf. [Se]).
\medskip
\flushpar
{\bf Theorem 1.2.1.} {\it A link group $G_{L}$ is $p$-good.}
\medskip
\flushpar
{\it Proof.} We shall say that a subgroup $G$ of finite index in $G_{L}$
 is \text{\it open} if $[G_{L}:G]$ is a power of $p$.
Let $M$ be a finite $p$-primary $\widehat{G_{L}}$-module.
We shall show by induction on the length of $M$ that if
$G$ is an open subgroup of $G_{L}$ then there is a smaller open subgroup
$G_1$ such that restriction from $H^2(G,M)$ to $H^2(G_1,M)$ is trivial.

Suppose first that $M={\bold F}_{p}$, with trivial $G_{L}$-action,
and let $H^*(G)$ denote $H^*(G,{\bold F}_{p})$, for ease of reading.
Since $[G_{L}:G]$ is finite and $G_{L}/G_{L}^{(2)}\cong {\bold Z}^{n}$
there is an epimorphism $\tau:G\to C= {\bold Z}/p{\bold Z}$.
Then $K=\roman{Ker}(\tau)$ is another open subgroup of $G_{L}$.
The Hochschild-Serre spectral sequence for $G$ as an extension of $C$ by
$K$ has $E_2$ term $E_2^{p,q}=H^p(C,H^q(K))$, $r^{th}$ differential $d_r$ of
bidegree $(r,1-r)$ and converges to $H^*(G)$.
Since $H^{p}(K)=0$ for $p>2$ there are only three nonzero rows,
and since $H^*(G)=0$ for $*>2$ we see that $d^{p,2}_3$ is an isomorphism
for all $p\geq1$. The spectral sequence is an algebra over the ring
$H^*(C)=E_2^{*,0}$.Since $C$ has cohomological period 2, the cup product 
with a generator of $H^2(C)\cong{\bold F}_p$ induces isomorphisms
$\gamma_{2}^{p,q}:E_2^{p,q}\cong{E}_2^{p+2,q}$ such that
$d_2^{p+2,q}\gamma_{2}^{p,q}=\gamma_{2}^{p+2,q-1}d_2^{p,q}$ for all $p,q\geq0$. 
Therefore we have the isomorphisms $\text{Ker}(d_{2}^{p,q}) \simeq
\text{Ker}(d_{2}^{p+2,q}), \; \text{Im}(d_{2}^{p,q}) \simeq
\text{Im}(d_{2}^{p+2,q}) \; \text{for any} \, p ,q\geq 0$. In particular
we have the isomorphisms $\gamma_{3}^{0,2} : E_{3}^{0,2} =
\text{Ker}(d_{2}^{0,2}) \cong  E_{3}^{2,2} = \text{Ker}(d_{2}^{2,2})$
and $\gamma_{3}^{3,0} : E_{3}^{3,0} =
E_{2}^{3,0}/\text{Im}(d_{2}^{1,1}) \cong E_{3}^{5,0} =
E_{2}^{5,0}/\text{Im}(d_{2}^{3,1})$ with $d_{3}^{2,2} \circ
\gamma_{3}^{0,2} = \gamma_{3}^{2,1}\circ d_{3}^{0,2}$. It follows that
$d^{0,2}_3$ is also an isomorphism, and so $E_\infty^{0,2}=0$.
But the edge homomorphism from $H^2(G)$ to $H^2(K)$ factors through
$E^{0,2}_\infty\leq E^{0,2}_2=H^2(K)^C$, and so is 0.

In general, $M$ has a finite composition series whose factors are copies
of the simple module ${\bold F}_{p}$.
Suppose that $M_1$ is a maximal proper submodule of $M$,
with quotient $M/M_1\cong {\bold F}_{p}$.
Restriction from $G$ to $K$ induces a homomorphism from the exact
sequences of cohomology corresponding to the coefficient sequence
$0\to M_1\to M\to {\bold F}_{p}\to0$. The result for ${\bold F}_{p}$ 
implies that the image of $H^2(G;M)$ lies in the image of $H^2(K,M_1)$.
By the hypothesis of induction we may assume the result is true for
$M_1$, and so there is an open subgroup $K_1<K$ such that restriction 
from $H^2(K,M_1)$ to $H^2(K_1,M_1)$ is trivial.
Hence restriction from $H^2(G,M)$ to $H^2(K_1,M)$ is also trivial.
This establishes the inductive step.

In particular, restriction from $H^2(G_{L},M)$ to $H^2(J,M)$ is trivial,
for some open subgroup $J$, and so the result follows,
as in Exercise 1 of Chapter I.\S2.6 of [Se].
(This exercise is stated in terms of profinite completions,
but extends easily to the pro-$p$ case).
\qed
\medskip
\flushpar
Since the cohomological dimension $cd (G_{L}) \leq 2$, with equality if
and only if $L$ is nontrivial, Theorem 1.2.4 gives the corresponding
bound for the pro-$p$ completion $\widehat{G_{L}}$.
\medskip
\flushpar
{\bf Corollary 1.2.2.} {\it The cohomological $p$-dimension
$cd_{p}(\widehat{G_{L}}) \leq 2$.}
\medskip
\flushpar
If the Milnor invariants of $L$ are all 0 mod $p$ (cf. Section 2), then
$\widehat{G_{L}}$  is a free pro-$p$ group and so
$cd_p(\widehat{G_{L}})$ may be strictly less than $cd(G_{L})$. In
particular, this is so if $L$ is a nontrivial knot.
\medskip
\flushpar
\head {\bf 2. $p$-adic Milnor invariants} \endhead
\medskip
\flushpar
2.1.{\it The pro-$p$ Fox differential calculus.} Let $\widehat{F}$ be
the pro-$p$ completion of the free group $F$ on $n$ generators
$x_{1},\cdots ,x_{n}$. Y. Ihara [Ih] extended the Fox differential
calculus on the abstract free group $F$ ([F]) to that on $\widehat{F}$.
The basic result is stated as the following
\medskip
\flushpar
{\bf Theorem 2.1.1} ([Ih]).  {\it There exists a unique continuous
${\bold Z}_{p}$-homomorphism }
$$  \partial_{i} = \frac{\partial}{\partial x_{i}} \; : \; {\bold
Z}_{p}[[\widehat{F}]] \longrightarrow {\bold Z}_{p}[[\widehat{F}]]$$
{\it for each $i$ ($1 \leq i \leq n$) such that any element $\alpha \in
{\bold Z}_{p}[[\widehat{F}]]$ is expressed uniquely in the form }
$$ \alpha = \epsilon(\alpha)1 + \sum_{i=1}^{n}\partial_{i}(\alpha)(x_{i}
- 1)$$
{\it where $\epsilon$ is the augmentation map ${\bold
Z}_{p}[[\widehat{F}]] \rightarrow {\bold Z}_{p}.$}
\medskip
\flushpar
The higher order derivatives are defined inductively by
$$ \partial_{i_{1}}\cdots \partial_{i_{r}}(\alpha) =
\partial_{i_{1}}(\partial_{i_{2}}\cdots \partial_{i_{r}}(\alpha)).$$
Here are some basic rules (cf. [Ih,2]).
\smallskip
\flushpar
(2.1.2) 1. If one restricts $\partial_{i}$ to  ${\bold Z}[F]$ under the
natural embedding ${\bold Z}[F] \rightarrow {\bold
Z}_{p}[[\widehat{F}]]$, we get the usual Fox derivative on ${\bold
Z}[F]$ ([F]).
\flushpar
2. $\partial_{i}(x_{j}) = \delta_{ij}$ (Kronecker delta).
\flushpar
3. $ \partial_{i}(\alpha \beta) = \partial_{i}(\alpha)\epsilon(\beta) +
\alpha \partial_{i}(\beta)\; $   ($\alpha, \beta \in  {\bold
Z}_{p}[[\widehat{F}]]$).
\flushpar
4. $\partial_{i}(f^{-1}) =  - f^{-1}\partial_{i}(f)\;$    ($f \in
\widehat{F}$).
\flushpar
5. For $f \in \widehat{F}$ and $a \in {\bold Z}_{p}$,
$\partial_{i}(f^{a}) = b \partial_{i}(f)$, where $b$ is any element of
${\bold Z}_{p}[[\widehat{F}]]$ such that $ b (f - 1) = f^{a} - 1$.
\flushpar
6. Let $\widehat{F'}$ be another free pro-$p$ group on  $x_{1}',\cdots
,x_{m}'$ and let $\varphi : \widehat{F} \rightarrow \widehat{F'}$ be a
continuous surjective homomorphism. Then one has
$\partial_{i}^{'}(\varphi(\alpha)) = \sum_{j=1}^{n}
\varphi(\partial_{j}(\alpha))\partial_{i}^{'}(\varphi(x_{j}))$, where
$\partial_{i}^{'} = \frac{\partial}{\partial x_{i}'}, \alpha \in {\bold
Z}_{p}[[\widehat{F}]].$
\flushpar
7. For $f \in \widehat{F}$, $\epsilon( \partial_{i}^{r}(f)) = {a
\choose r}$ where $a = \epsilon(\partial_{i}(f))$ and ${a \choose r} =
\frac{a(a-1)\cdots (a-r+1)}{r !} \in {\bold Z}_{p}$.
\medskip

 Let ${\bold Z}_{p}\langle \langle X_{1},\cdots ,X_{n}\rangle \rangle$
 be the formal power series ring over ${\bold Z}_{p}$ in non-commuting
variables $X_{1},\cdots ,X_{n}$ which is compact in the topology taking
the ideals $I(r)$  of power series with homogeneous components of degree
$\geq r$ as the system of neighborhood of $0$. The pro-$p$ Magnus
embedding $M$ of $\widehat{F}$ into ${\bold Z}_{p}\langle\langle
X_{1},\cdots ,X_{n}\rangle \rangle^{\times}$ is defined by
$$ M(x_{i}) = 1+X_{i}, \; M(x_{i}^{-1}) = 1-X_{i}+X_{i}^{2}+\cdots $$
and it is extended to give the isomorphism  ${\bold
Z}_{p}[[\widehat{F}]] \simeq  {\bold Z}_{p}\langle\langle X_{1},\cdots
,X_{n}\rangle \rangle$ of compact ${\bold Z}_{p}$-algebras. The
resulting expansion of $\alpha \in {\bold Z}_{p}[[\widehat{F}]]$ is
given by the Fox derivatives:
$$  \split M(\alpha)  =  & \; \epsilon(\alpha) + \sum_{I=(i_{1}\cdots
i_{r})}\epsilon_{I}(\alpha)X_{i_{1}}\cdots X_{i_{r}}, \\
&  \epsilon_{I}(\alpha) = \epsilon(\partial_{i_{1}}\cdots
\partial_{i_{r}}(\alpha)), \;\; I = (i_{1}\cdots i_{r}).  \endsplit
 \tag 2.1.3 $$
Finally, we recall the following fact ([Ko],7.14): Let
$\widehat{J_{q}}$ be the two-sided ideal of ${\bold
Z}_{p}[[\widehat{F}]]$ generated by $q$ and the augmentation ideal
$I_{{\bold Z}_{p}[[\widehat{F}]]}$. Then for $f \in \widehat{F}$ and $k
\geq 1$, we have
$$ \align f \in \widehat{F}^{(k,q)} & \Longleftrightarrow  f -1 \in
\widehat{J_{q}} \tag 2.1.4\\
        & \Longleftrightarrow M(f) = 1 + (\text{term of degree} \geq k)
\;\; (\text{for case} \; q = 0).\\
\endalign $$
\medskip
\flushpar
2.2. {\it $p$-adic Milnor invariants.} We keep the same notation as in
Section 1. Let $\widehat{F}$ be the free pro-$p$ group on $x_{1},\cdots
,x_{n}$ where each $x_{i}$ represents the $i$-th meridian. For a
multi-index $I = (i_{1}\cdots i_{r}), r \geq 1$, we define the {\it
$p$-adic Milnor number} $\hat{\mu}(I)$ by
$$ \split \hat{\mu}(I) & := \epsilon_{I'}(y_{i_{r}}) \;\; I' =
(i_{1}\cdots i_{r-1})\\
                     & = \epsilon(\partial_{i_{1}}\cdots
\partial_{i_{r-1}}(y_{i_{r}}))
\endsplit
\tag 2.2.1 $$
where $y_{j} \in \widehat{F}$ represents the $j$-th ``longitude" in
$\widehat{G_{L}}$ (cf. Section 1.1). By convention, we set $\hat{\mu}(I)
= 0$ for $|I| = r =1$. We let $\hat{\Delta}(I)$ denote the ideal of
${\bold Z}_{p}$ generated by $\hat{\mu}(J)$ where $J$ runs over all
cyclic permutations of proper subsequences of $I$. We then define the
{\it $p$-adic Milnor invariant} $\overline{\hat{\mu}}(I)$ by
$$ \overline{\hat{\mu}}(I) := \hat{\mu}(I) \; \text{mod} \;
\hat{\Delta}(I). \tag 2.2.2$$
Since the usual Milnor number $\mu (I), I = (i_{1}\cdots i_{r})$ is
defined by $\epsilon(\partial_{i_{1}}\cdots
\partial_{i_{r-1}}(y_{i_{r}}^{(r)}))$ by (2.1.1),1,  Remark 1.1.3, (1)
and (2.2.1) yield
$$ \hat{\mu}(I) = \mu(I) \;\; \text{and} \;\; \hat{\Delta}(I) =
\Delta(I) \tag 2.2.3$$
as elements in ${\bold Z}_{p}$. Hence, we have
\smallskip
\flushpar
(2.2.4)  $\overline{\hat{\mu}}(I)$'s are isotopy invariants of $L$ and
satisfy the same properties such as the cycle symmetry and Shuffle
relations $\overline{\mu}(I)$'s enjoy ([Mi]).
\smallskip
\flushpar
Theorem 1.1.2 and (2.1.4) implies the following:
\smallskip
\flushpar
(2.2.5) All $\overline{\hat{\mu}}(I) = 0$ for $|I| < r$ if and only if
$\widehat{\pi} : \widehat{F} \rightarrow \widehat{G_{L}}$ induces an
isomorphism $\widehat{F}/\widehat{F}^{(r,q)} \simeq
\widehat{G_{L}}/\widehat{G_{L}}^{(r,q)}$. In particular, if all $p$-adic
Milnor invariants of $L$ are zero, one has $\widehat{G_{L}} \simeq
\widehat{F}$. This is the case for boundary links.
\smallskip
\flushpar
Finally, we remark that (2.1.2), 7 implies
\smallskip
\flushpar
(2.2.6)\hskip2cm
$\displaystyle{\hat{\mu}(\undersetbrace{\text{r}}\to{i\cdots i}j) =
{\text{lk}(K_{i},K_{j}) \choose r} \;\; \text{for} \; i\neq j.}$
\medskip

The Milnor invariant is also given by the Massey products in the
cohomology of $\widehat{G_{L}}$. For the normalized Massey system in
profinite group cohomology and sign convention, we refer to [Mo3]. Let
$\xi_{1},\cdots ,\xi_{n}$ be the ${\bold Z}_{p}$-basis of
$H^{1}(\widehat{G_{L}},{\bold Z}_{p})$ dual to the meridians $m_{i}$'s,
and let $\eta_{j} \in H_{2}(\widehat{G_{L}},{\bold Z}_{p})$ be the image
of $[x_{j},y_{j}]$ under the transgression $H_{1}(\widehat{N},{\bold
Z}_{p})_{\widehat{G_{L}}} \rightarrow H_{2}(\widehat{G_{L}},{\bold
Z}_{p})$. Then for $I = (i_{1}\cdots i_{r}), r\geq 2$, there is a
normalized Massey system $M$ for the product $\langle \xi_{i_{1}},\cdots
,\xi_{i_{r}} \rangle \in H_{2}(\widehat{G_{L}},{\bold
Z}_{p}/\hat{\Delta}(I))$ so that
$$
\langle \xi_{i_{1}},\cdots ,\xi_{i_{r}} \rangle(\eta_{j}) =
\cases
(-1)^{r}\overline{\hat{\mu}}(I) & \quad j = i_{r} \neq i_{1}, \\
(-1)^{r+1}\overline{\hat{\mu}}(I) & \quad j = i_{1} \neq i_{r},\\
0 & \quad \text{otherwise}.
\endcases
\tag 2.2.7
$$
\medskip
\flushpar
\head {\bf 3. Completed Alexander modules} \endhead
\medskip
\flushpar
3.1. {\it The Alexander module of $\widehat{G_{L}}$.} The Alexander
module of a finitely presented pro-$p$ group was introduced in several
modes in [Mo2]. As a particular case, the Alexander module of the
pro-$p$ link group $\widehat{G_{L}}$ is defined using the pro-$p$ Fox
free differential calculus as follows. We keep the same notation as in
Section 1 and 2. Let $\widehat{\psi}$ be the abelianization map
$\widehat{G_{L}} \rightarrow \widehat{H} :=
\widehat{G_{L}}/\widehat{G_{L}}^{(2)} = {\bold Z}_{p}^{n}$ and denote by
the same $\widehat{\psi}$ the continuous ${\bold Z}_{p}$-homomorphism
${\bold Z}_{p}[[\widehat{G_{L}}]] \rightarrow {\bold
Z}_{p}[[\widehat{H}]]$ on the completed group rings. We identify ${\bold
Z}_{p}[[\widehat{H}]]$ with the commutative formal power series ring
$\widehat{\Lambda_{n}} := {\bold Z}_{p}[[X_{1},\cdots ,X_{n}]]$ over
${\bold Z}_{p}$ by setting $t_{i} := \widehat{\psi}\circ \widehat{\pi}
(x_{i}) = 1+X_{i}$. By Theorem 1.2.1, we call the matrix over
$\widehat{\Lambda_{n}}$
$$ \widehat{P_{L}} = \left( \widehat{\psi}\circ
\widehat{\pi}(\partial_{j}([x_{i},y_{i}]))\right) \tag 3.1.1 $$
 the {\it Alexander matrix} of $\widehat{G_{L}}$. We then define the
{\it Alexander module} $\widehat{A_{L}}$ of $\widehat{G_{L}}$ by the
compact $\widehat{\Lambda_{n}}$-module presented by $\widehat{P_{L}}$:
$$ \widehat{A_{L}} := \text{Coker} (\widehat{\Lambda_{n}}^{n} \overset
\widehat{P_{L}} \to \longrightarrow \widehat{\Lambda_{n}}^{n} ) \tag
3.1.2$$
and call $\widehat{A_{L}}$ the {\it completed Alexander module} of $L$
over $\widehat{\Lambda_{n}}$. Let $\widehat{\Lambda} = {\bold
Z}_{p}[[X]]$, the {\it Iwasawa algebra} and $\tau :
\widehat{\Lambda_{n}} \rightarrow \widehat{\Lambda}$ the reducing
homomorphism defined by $\tau(X_{i}) = X$ $(1\leq i\leq n)$. Then the
{\it reduced completed Alexander module} $\widehat{A_{L}}^{red}$ is
defined by the compact $\widehat{\Lambda}$-module
$$ \widehat{A_{L}}^{red} := \tau(\widehat{A_{L}}) =
\widehat{A_{L}}\otimes_{\widehat{\Lambda_{n}}}\widehat{\Lambda} \tag
3.1.3$$
which is presented by $\tau(\widehat{P_{L}}) = \widehat{P_{L}}(X,\cdots
,X)$. For a knot $K$, $\widehat{P_{L}} = O_{n}$ (zero matrix) and so we
have
$$ \widehat{A_{L}} =  \widehat{A_{L}}^{red} = \widehat{\Lambda}. \tag
3.1.4$$
We also define the {\it $i$-th completed Alexander ideal}
$\widehat{E_{i}}(L)$ of $L$ by the $i$-th elementary ideal
$E_{i}(\widehat{A_{L}})$ of $\widehat{A_{L}}$  and {\it $i$-th $p$-adic
Alexander series} $\widehat{\Delta_{i}}(L)$ of $L$ by the greatest
common divisor $\Delta_{i}(\widehat{E_{i}}(L))$ of generators of the
ideal $\widehat{E_{i}}(L)$:
$$ \widehat{E_{i}}(L) := E_{i}(\widehat{A_{L}}), \;\;
\widehat{\Delta_{i}}(L) := \Delta_{i}(\widehat{E_{i}}(L)) \tag 3.1.5 $$
\smallskip

The relation with the usual Alexander module is given as follows. Let
$\psi : G_{L} \rightarrow H = H_{1}(X_{L},{\bold Z})$ be the
abelianization map which induces the ring homomorphism $\psi : {\bold
Z}[G_{L}] \rightarrow {\bold Z}[H]$ on the group rings where ${\bold
Z}[H]$ is identified with the Laurent polynomial ring $\Lambda_{n} =
{\bold Z}[t_{1}^{\pm 1},\cdots ,t_{n}^{\pm 1}]$, $t_{i} = \psi \circ \pi
(x_{i})$. Given a presentation $ G_{L} = \langle x_{1},\cdots ,x_{m}\;
|\; r_{1}=\cdots r_{l}=1 \rangle$ $(m \geq n)$, the Alexander module of
$L$ is given as the $\Lambda_{n}$-module presented by the Alexander
matrix $\left(\psi\circ \pi (\partial_{j}(r_{i}))\right)$. By Theorem
1.1.1, we can take $m = n$ and the relators to be
$[x_{i},y_{i}^{(k+1)}]$ ($1\leq i \leq n$) and some finite number of
generators $f_{i}^{(k+1)}$ of $F^{(k+1,p)}$ $(k\geq 1)$. We let $J_{p} =
\widehat{J_{p}} \cap {\bold Z}[F]$ and $J_{\Lambda_{n},p} = \psi\circ
\pi (J_{p})$. Since $\partial_{j}(f_{i}^{(k+1)}) \in {J_{p}}^{k}$ by
(2.1.4), passing to the quotients modulo $(J_{\Lambda_{n},p})^{k}$,
$A_{L}/(J_{\Lambda_{n},p})^{k}A_{L} =
A_{L}\otimes_{\Lambda_{n}}\Lambda_{n}/(J_{\Lambda_{n},p})^{k}$ is
presented by the matrix $\left( \psi \circ \pi
(\partial_{j}([x_{i},y_{i}^{(k+1)}]))  \; \text{mod} \;
(J_{\Lambda_{n},p})^{k} \right)$. Here we see by  Remark 1.1.3, (1) that
the elements $\{ \psi \circ \pi (\partial_{j}([x_{i},y_{i}^{(k+1)}])) \;
\text{mod} \;  (J_{\Lambda_{n},p})^{k} \}$ form an inverse system with
respect to $k$ and its limit is given as $\widehat{\psi}\circ
\widehat{\pi}(\partial_{j}([x_{i},y_{i}]))$  under the identification
$\varprojlim\limits_{k} \Lambda_{n}/(J_{\Lambda_{n},p})^{k} =
\widehat{\Lambda_{n}}$ defined by $t_{i} = 1+X_{i}$. Hence by (3.1.2) we
have
$$ \widehat{A_{L}} = \varprojlim\limits_{k} A_{L}\otimes_{\Lambda_{n}}
(\Lambda_{n}/(J_{\Lambda_{n,p}})^{k}) =
A_{L}\otimes_{\Lambda_{n}}\widehat{\Lambda_{n}}. \tag 3.1.6 $$
Similarly, $\widehat{A_{L}}^{red}$ is related with the usual reduced
Alexander module $A_{L}^{red} = \tau(A_{L})$ by
$$ \widehat{A_{L}}^{red} = A_{L}\otimes_{\Lambda}\widehat{\Lambda} \tag
3.1.7$$
where $\Lambda = {\bold Z}[t^{\pm 1}]$ is embedded into
$\widehat{\Lambda}$ by $t = 1+X$.
\medskip
\flushpar
3.2. {\it The $p$-adic Traldi  matrix.} The Alexander matrix
$\widehat{P_{L}}$ of $\widehat{G_{L}}$ (3.1.1) is computed explicitly as
a {\it universal $p$-adic higher linking matrix} in terms of $p$-adic
Milnor numbers. This is regarded  as a $p$-adic strengthening of the
work by L. Traldi [Tr].
\medskip
\flushpar
{\bf Definition 3.2.1.} The {\it $p$-adic Traldi matrix}
$\widehat{T_{L}} = (\widehat{T_{L}}(i,j))$ of $L$ over
$\widehat{\Lambda_{n}}$ is defined by
$$ \widehat{T_{L}}(i,j) = \cases - \displaystyle{\sum_
{r \geq 1}\sum  \Sb 1\leq i_{1},\cdots ,i_{r}\leq n \\ i_{r}\neq i
\endSb
\hat{\mu}(i_{1}\cdots i_{r}i)X_{i_{1}}\cdots X_{i_{r}}} & \quad i=j \\
\displaystyle{\hat{\mu}(ji)X_{i}+\sum_{r \geq 1}\sum_{1\leq i_{1}\cdots
i_{r}\leq n}\hat{\mu}(i_{1}\cdots i_{r}ji)X_{i}X_{i_{1}}\cdots
X_{i_{r}}} & \quad i\neq j.
\endcases $$
and  we also define the {\it reduced $p$-adic Traldi matrix}
$\widehat{T_{L}}^{red}$ of $L$ over $\widehat{\Lambda}$ by
$$\widehat{T_{L}}^{red} := \tau(\widehat{T_{L}}) =
\widehat{T_{L}}(X,\cdots ,X). $$
Our theorem is then stated as
\medskip
\flushpar
{\bf Theorem 3.2.2.} {\it The $p$-adic Traldi matrix $\widehat{T_{L}}$
 gives a presentation matrix for the completed Alexander module
$\widehat{A_{L}}$  over $\widehat{\Lambda_{n}}$, and the reduced
$p$-adic Traldi matrix $\widehat{T_{L}}^{red}$ gives a presentation
matrix for the reduced completed Alexander module
$\widehat{A_{L}}^{red}$ over $\widehat{\Lambda}$.}
\medskip
\flushpar
{\it Proof.} By (3.1.1),(3.1.2) and (3.1.3), it suffices to show
$\widehat{\psi}\circ \widehat{\pi}(\partial_{j}([x_{i},y_{i}]))
= \widehat{T_{L}}(i,j)$. By the rules $2\sim 4$ of (2.1.2), we have
$$ \partial_{j}([x_{i},y_{i}])  = (1-x_{i}y_{i}x_{i}^{-1})\delta_{ij} +
x_{i}(1-y_{i}x_{i}^{-1}y_{i}^{-1})\partial_{j}( y_{i}). $$
Here  $y_{i} = 1+\displaystyle{\sum_{r\geq 1}\sum_{1\leq i_{1},\cdots
,i_{r}\leq n}\hat{\mu}(i_{1}\cdots i_{r}i)(x_{i_{1}}-1)\cdots
(x_{i_{r}}-1)}$ by (2.1.3) and (2.2.1). Hence we get
$$  \align \widehat{\psi}\circ\widehat{\pi}
\left(\partial_{j}([x_{i},y_{i}])\right) &
= \delta_{i,j} - \sum_{r\geq 1}\sum_{1\leq i_{1},\cdots ,i_{r}\leq
n}\hat{\mu}(i_{1}\cdots i_{r}i)X_{i_{1}}\cdots X_{i_{r}}  \tag 3 \\
 & +\hat{\mu}(ji)X_{i}+ \sum_{r\geq 1}\sum_{1\leq i_{1},\cdots
,i_{r}\leq n}\hat{\mu}(i_{1}\cdots i_{r}ji)X_{i}X_{i_{1}}\cdots
X_{i_{r}}. \endalign  $$
which yields the assertion. $\qed$
\medskip
\flushpar
The following is an extension of (3.1.4).
\medskip
\flushpar
{\bf Corollary 3.2.3.} {\it For a link whose $p$-adic Milnor invariants
are all zero, we have}
$$ \widehat{A_{L}} \simeq {\widehat{\Lambda_{n}}}^{n}, \;\;
\widehat{A_{L}}^{red} \simeq {\widehat{\Lambda}}^{n}.$$
{\it This is the case for boundary links.}
\medskip

Finally,  we introduce the {\it truncated $p$-adic Traldi matrices.}
\medskip
\flushpar
{\bf Definition 3.2.4.} For $k \geq 2$, the {\it $k$-th truncated
$p$-adic Traldi matrix}
$\widehat{T_{L}}^{(k)} = (\widehat{T_{L}}^{(k)}(i,j))$  is defined by
$$ \widehat{T_{L}}^{(k)}(i,j) =
\cases -\displaystyle{\sum_{r=1}^{k-1} \sum \Sb 1\leq i_{1},\cdots
,i_{r}\leq n \\ i_{r}\neq i \endSb \hat{\mu}(i_{1}\cdots
i_{r}i)X_{i_{1}}\cdots X_{i_{r}}} &  \quad i = j \\
 \hat{\mu}(ji)X_{i} + \displaystyle{\sum_{r=1}^{k-2}\sum_{1\leq
i_{1},\cdots ,i_{r}\leq n} \hat{\mu}(i_{1}\cdots
i_{r}ji)X_{i}X_{i_{1}}\cdots X_{i_{r}}} & \quad i \neq j
\endcases $$
and we also define the {\it $k$-th truncated reduced $p$-adic Traldi
matrix} $\widehat{T_{L}}^{red, (k)}$  by
$$ \widehat{T_{L}}^{red, (k)} := \tau(\widehat{T_{L}}^{(k)}) =
\widehat{T_{L}}^{(k)}(X,\cdots ,X). $$
\medskip
\flushpar
We note that $\widehat{T_{L}}^{red,(2)}$ is the {\it linking matrix}
multiplied by $X$, where the linking matrix $C = (C_{ij})$ is defined by
$C_{ii} = -\sum_{j\neq i}\text{lk}(K_{i},K_{j})$ and $C_{ij} =
\text{lk}(K_{i},K_{j})$ for $i\neq j$. Thus the $p$-adic Traldi matrix
$\widehat{T_{L}}$ is regarded as a {\it universal higher  linking
matrix} over $\widehat{\Lambda_{n}}$ which contains all information on
the completed Alexander module. In the following section, we derive from
$\widehat{T_{L}}$ the information on the $p$-homology groups of
$p^{m}$-fold cyclic branched covers along $L$.
\bigskip
\flushpar
\head{\bf 4. Galois module structure for the $p$-homology group of a
$p$-fold cyclic branched cover} \endhead
\medskip
\flushpar
4.1. {\it Galois module structure of the $p$-homology of a $p$-fold
cover.} Let $X_{\infty}$ be the infinite cyclic cover of $X_{L}$
 associated to the kernel of the homomorphism $G_{L} \rightarrow \langle
t \rangle$ sending each meridian to $t$. Let $M$ be the the completion
of the $p$-fold subcover of $X_{\infty}$ over $X_{L}$ so that $M$ is a
$p$-fold cyclic cover of $S^{3}$ branched along $L$.  We set $\nu_{d}(t)
= t^{d-1}+\cdots + t + 1$ for $d\geq 1$. Let $\phi : M \rightarrow
S^{3}$ be the covering map and $\sigma$  a generator of its Galois
group. Since $\nu_{p}(\sigma) = tr\circ \phi_{*}=0$ where $tr :
H_{1}(S^{3},{\bold Z}) \rightarrow H_{1}(M,{\bold Z})$ is the transfer,
we can regard $H_{1}(M,{\bold Z})$ as a module over the Dedekind ring
${\Cal O} = {\bold Z}[\langle \sigma \rangle ]/(\nu_{p}(\sigma)) =
{\bold Z}[\zeta]$, $\zeta = \sigma \; \text{mod} \; (\nu_{p}(\sigma))$
 is a primitive $p$-th root of 1. Hence $H_{1}(M,{\bold Z}_{p}) =
H_{1}(M,{\bold Z})\otimes_{\bold Z}{\bold Z}_{p}$ is regarded as a
module over the complete discrete valuation ring $\widehat{\Cal O} =
{\Cal O}\otimes_{\bold Z}{\bold Z}_{p} = {\bold Z}_{p}[\zeta]$. Note
that $\widehat{\Cal O}$ is the completion of ${\Cal O}$ with respect to
the maximal ideal $\frak{p}$ generated by the prime element $\pi :=
\zeta -1$ and the residue field $\widehat{\Cal O}/\frak{p}$ is ${\bold
F}_{p}$. By Theorem 2.2.1, we can derive the following information on a
presentation matrix for the $\widehat{\Cal O}$-module $H_{1}(M,{\bold
Z}_{p})$. Note that the evaluation of a power series $F(X) \in {\bold
Z}_{p}[[X]]$ at $ s = \pi$ makes sense in the $\frak{p} =
(\pi)$-adically complete ring $\widehat{\Cal O}$.
\medskip
\flushpar
{\bf Theorem 4.1.1.} {\it A presentation matrix for $H_{1}(M,{\bold
Z}_{p})\oplus \widehat{\Cal O}$ over $\widehat{\Cal O}$ is given by
$\widehat{T_{L}}^{red}(\pi)$.  Further, for any integer $k \geq 2$, a
presentation matrix for $(H_{1}(M,{\bold Z}_{p})\otimes_{\widehat{\Cal
O}} \widehat{\Cal O}/\frak{p}^{k})\oplus \widehat{\Cal O}/\frak{p}^{k}$
 over $\widehat{\Cal O}/\frak{p}^{k}$ is given by $\widehat{T_{L}}^{red,
(k)}(\pi)$. Here  $\widehat{T_{L}}^{red}$ (resp.
$\widehat{T_{L}}^{red, (k)}$) is the reduced (resp. reduced truncated)
Traldi matrix  defined in Section 3.2.}
\medskip
\flushpar
{\it Proof.} Note that the well-known relation ([S1,Theorem 6])
$$ H_{1}(M,{\bold Z}) \simeq H_{1}(X_{\infty},{\bold
Z})/\nu_{p}(t)H_{1}(X_{\infty},{\bold Z})$$
 is an $\widehat{\Cal O}$-isomorphism since $\sigma$ acts on the r.h.s
by $t$. Hence we have the following isomorphisms over $\widehat{\Cal O}$

$$ \align
H_{1}(M,{\bold Z}_{p}) & \simeq (H_{1}(X_{\infty},{\bold
Z})/\nu_{p}(t)H_{1}(X_{\infty},{\bold Z}))\otimes_{\bold Z}{\bold Z}_{p}
\\
& \simeq H_{1}(X_{\infty},{\bold
Z})\otimes_{\Lambda}(\Lambda\otimes_{\bold Z}{\bold
Z}_{p})/(\nu_{p}(t))) \\
& \simeq  H_{1}(X_{\infty},{\bold
Z})\otimes_{\Lambda}(\widehat{\Lambda}/(\nu_{p}(1+X))) \\
& \simeq H_{1}(X_{\infty},{\bold Z})\otimes_{\Lambda}\widehat{\Cal O}.
\endalign$$
Since $A_{L}^{red} \simeq H_{1}(X_{\infty},{\bold Z})\oplus \Lambda$ as
$\Lambda$-module ([H, 5.4]), tensoring with $\widehat{\Cal O}$ over
$\Lambda$, we have an isomorphism of $\widehat{\Cal O}$-modules
$$ A_{L}^{red}\otimes_{\Lambda}\widehat{\Cal O} \simeq H_{1}(M,{\bold
Z}_{p})\oplus\widehat{\Cal O}.$$
Since the l.h.s is same as
$\widehat{A_{L}}^{red}\otimes_{\widehat{\Lambda}}\widehat{\Cal O}$, the
first assertion follows from Theorem 3.2.4.  The second assertion is
obtained from the first one by taking modulo $\frak{p}^{k}$. $\qed$
\medskip
\flushpar
Now, we asuume that $H_{1}(M,{\bold Z})$ is finite so that
$H_{1}(M,{\bold Z}_{p})$ is the $p$-primary part  of $H_{1}(M,{\bold
Z})$. Using Theorem 3.1.1, we will see the $\widehat{\Cal O}$-module
structure of $H_{1}(M,{\bold Z}_{p})$ more precisely. First, we recall
the following result on the $\frak{p}$-rank of $H_{1}(M,{\bold Z}_{p})$
 (cf. [Mo1],[Rez]).
\medskip
\flushpar
{\bf Lemma 4.1.2.} {\it $H_{1}(M,{\bold Z}_{p})\otimes_{\widehat{\Cal
O}}{\bold F}_{p}$ has dimension $n-1$ over ${\bold F}_{p}$.}
\medskip
\flushpar
{\it Proof.} By [Mo1], the map $\Phi : H_{1}(M,{\bold Z}) \rightarrow
{\bold F}_{p}^{n}$ defined by
$\Phi(c) := (\text{lk}(\phi_{*}(c),K_{i}) \; \text{mod} \; p)$ induces
an isomorphism
$$ H_{1}(M,{\bold Z})/(\sigma -1)H_{1}(M,{\bold Z}) \simeq \{ (\xi_{i})
\in {\bold F}_{p}^{n} \; | \; \sum_{i=1}^{n}\xi_{i} = 0 \} \simeq {\bold
F}_{p}^{n-1}$$
where the l.h.s is $H_{1}(M,{\bold Z})\otimes_{\Cal O}{\Cal O}/{\frak p}
= H_{1}(M,{\bold Z}_{p})\otimes_{\widehat{\Cal O}}{\bold F}_{p},$ and
hence we are done. $\qed$
\medskip
\flushpar
By Lemma 4.1.2, $H_{1}(M,{\bold Z}_{p})$ has the form
$$ H_{1}(M,{\bold Z}_{p}) = \bigoplus_{i=1}^{n-1}\widehat{\Cal O}/{\frak
p}^{a_{i}} \;\;\; \; \; (a_{i} \geq 1)$$
as $\hat{\Cal O}$-module. Hence, the determination of $\widehat{\Cal
O}$-module structure of $H_{1}(M,{\bold Z}_{p})$ is equivalent to the
determination of the  $\frak{p}^{k}$-rank
$$ e_{k} := \# \{ i \; | \; a_{i} \geq k \}  = \text{dim}_{{\bold
F}_{p}}H_{1}(M,{\bold F}_{p})\otimes_{\widehat{\Cal
O}}(\frak{p}^{k-1}/\frak{p}^{k}) \;\; (k\geq 1).$$
We describe the ${\frak p}^{k}$-rank $e_{k}$ in terms of
$\widehat{T_{L}}^{red,(k)}(\pi)$. For an $n$ by $n$ matrix $A$ over
$\widehat{\Cal O}$, we denote by $A\otimes_{\widehat{\Cal
O}}(\frak{p}^{k-1}/\frak{p}^{k})$ the ${\bold F}_{p}$-linear map on the
${\bold F}_{p}$-vector space $(\frak{p}^{k-1}/\frak{p}^{k})^{n}$ induced
by $A$. Then Theorem 4.1.1 is restated as
\medskip
\flushpar
{\bf Theorem 4.1.3.} {\it Notation and assumption being as above, we
have}
$$ e_{k} = n - 1 - \text{rank}_{{\bold F}_{p}}(\widehat{T_{L}}^{red,
(k)}(\pi)\otimes_{\widehat{\Cal O}}(\frak{p}^{k-1}/\frak{p}^{k})),
\;\; k \geq 2.$$
Here we may call $\widehat{T_{L}}^{red, (k)}(\pi)$ the {\it $k$-th
higher linking matrix} in view of the following
\medskip
\flushpar
{\bf Corollary 4.1.4.} {\it For $k = 2$,  we have}
$$ e_{2} = n-1 - \; \text{rank}_{{\bold F}_{p}}(C \; \text{mod} \; p) $$

{\it where} $C = (C_{ij})$ {\it is the linking matrix of} $L$ {\it
defined by} $C_{ii} = -\sum_{j\neq i}\text{lk}(K_{i},K_{j})$ {\it and}
$C_{ij} = \text{lk}(K_{i},K_{j})$ {\it for} $i\neq j$.
\medskip
\flushpar
{\it Proof.} In fact, we have, by definition,
$$ \widehat{T_{L}}^{red,(2)}(\pi) = \pi C$$
and hence $\text{rank}_{{\bold
F}_{p}}(\widehat{T_{L}}^{red,(2)}\otimes_{\widehat{\Cal
O}}\frak{p}/\frak{p}^{2}) =  \text{rank}_{{\bold F}_{p}}(C \; \text{mod}
\; p).$  $\qed$
\medskip
\flushpar
4.2. {\it $2$-component case.} We suppose $n = 2$ and keep to assume
$H_{1}(M,{\bold Z})$ is finite. By Lemma 4.1.2,  $H_{1}(M,{\bold
Z}_{p})$ has the ${\frak p}$-rank $1$ so that we have
$$ H_{1}(M,{\bold Z}_{p}) = \widehat{\Cal O}/\frak{p}^{a}, \;\; a\geq
1.$$
 Hence $e_{k} = 0 \; \text{or}\; 1$ for $k\geq 2$, and by Theorem 4.1.3
we have
$$
e_{k} = 1  \Longleftrightarrow \widehat{T_{L}}^{red,(k)}(\pi) \equiv
O_{2} \; \text{mod}\; \pi^{k}.
$$
 Since $\widehat{T_{L}}^{red,(k)}(1,2)(\pi) =
-\widehat{T_{L}}^{red,(k)}(1,1)(\pi)$,
$\widehat{T_{L}}^{red,(k)}(2,2)(\pi) =
-\widehat{T_{L}}^{red,(k)}(2,1)(\pi)$, we have the following
\medskip
\flushpar
{\bf Theorem 4.2.1.}  {\it Suppose $n = 2$. For each $k\geq 1$, assuming
$e_{k} = 1$, we have}
$$
\align
e_{k+1} = 1 & \Longleftrightarrow
\cases
\displaystyle{\sum_{r=1}^{k}\sum_{i_{1},\cdots ,i_{r-1} =
1,2}\hat{\mu}(i_{1}\cdots i_{r-1}21)\pi^{r} \; \equiv \; 0 \;
\text{mod}\; \pi^{k+1}},\\
\displaystyle{\sum_{r=1}^{k}\sum_{i_{1},\cdots ,i_{r-1} =
1,2}\hat{\mu}(i_{1}\cdots i_{r-1}12)\pi^{r} \; \equiv \; 0 \;
\text{mod}\; \pi^{k+1}}. \\
\endcases
\endalign
$$
We give the condition in Theorem 4.2.1 in more concise forms for lower
$k$. In the following computation, we use simply the usual Milnor number
$\mu(I)$ instead of $\hat{\mu}(I)$ by (2.2.3).
\medskip
\flushpar
{\bf Example 4.2.2.} $e_{2}$: Since $e_{1} = 1$, we have by Theorem
4.2.1
$$ e_{2} = 1 \Longleftrightarrow \mu(12)\pi \equiv 0 \; \text{mod}\;
\pi^{2} \Longleftrightarrow \text{lk}(K_{1},K_{2}) \equiv 0 \;
\text{mod} \; p. \tag 4.2.2.1 $$
\medskip
\flushpar
$e_{3}$: Assume $\text{lk}(K_{1},K_{2}) \equiv 0$ mod $p$. By Theorem
4.2.1, we have
$$ \align
e_{3} = 1 & \Longleftrightarrow
\cases
 \mu(21)\pi + (\mu(121) + \mu(221))\pi^{2}  \; \equiv \; 0 \;
\text{mod}\; \pi^{3}, \\
\mu(12)\pi + (\mu(112) + \mu(212))\pi^{2}  \; \equiv \; 0 \;
\text{mod}\; \pi^{3}.
\endcases
\endalign
$$
By cycle symmetry, $ \mu(121) \equiv \mu(112), \; \mu(221) \equiv
\mu(212) \; \text{mod}\; \mu(12)$. Here  $\mu(112) \equiv \mu(221)
\equiv {\mu(12)\choose 2}$ mod $\mu(12)$ by (2.2.6). Thus we have
$\mu(121) + \mu(221) \equiv \mu(112) + \mu(212) \equiv \mu(12) \equiv 0$
 mod $p$. Hence, we have
$$ e_{3} = 1 \Longleftrightarrow \text{lk}(K_{1},K_{2})  \equiv  0 \;
\text{mod}\; p^{2}. \tag 4.2.2.2 $$
As one easily see, this condition is also equivalent to
$$ \mu(112) \equiv \mu(221) \equiv 0 \; \text{mod} \;  p. \tag 4.2.2.3
$$
\medskip
\flushpar
$e_{4}$: Assume $\text{lk}(K_{1},K_{2}) \equiv 0$ mod $p^{2}$. By
Theorem 4.2.1, we have
$$ \align
e_{4} = 1 & \Longleftrightarrow
\cases
\mu(21)\pi + (\mu(121)+\mu(221))\pi^{2} \\
\;\;\;\;\; + (\mu(1121)+\mu(1221)+\mu(2121)+\mu(2221))\pi^{3} \;
\equiv \; 0 \; \text{mod}\; \pi^{4}, \\
\mu(12)\pi + (\mu(112) + \mu(212))\pi^{2} \\
 \;\;\;\;\; + (\mu(1112)+\mu(1212)+\mu(2112)+\mu(2212))\pi^{3} \;
\equiv \; 0 \; \text{mod}\; \pi^{4}.
\endcases
\endalign
$$
As in case of $e_{3}$, we have $\mu(121) + \mu(221) \equiv \mu(12)
\equiv 0$ mod $p^{2}$. Similarly, $\mu(1121) \equiv \mu(1112) \equiv
{\mu(12)\choose 3}$ mod $\Delta(1121)$ and  $\mu(2221) \equiv
{\mu(12)\choose 3}$ mod $\Delta(2221)$ by (2.2.6). Since $\Delta(1121)
\equiv \Delta(2221) \equiv 0$ mod $p$, $\mu(1121)+\mu(2221) \equiv
\mu(12)/3 \equiv 0$ mod $p$. Finally, by shuffle relation, $\mu(1221) +
\mu(2121) + \mu(2211) \equiv \mu(2121) + 2\mu(2211) \equiv 0$ mod $p$.
Thus the first condition is equivalent to  $\mu(21)\pi -
\mu(2211)\pi^{3} \equiv 0$ mod $\pi^{4}$. Similarly, we see that
the second condition is equivalent to  $\mu(12)\pi -
\mu(1122)\pi^{3} \equiv 0$ mod $\pi^{4}$ which is same as the
first one. Hence, we obtain
$$ e_{4} = 1 \Longleftrightarrow \text{lk}(K_{1},K_{2}) -
\mu(1122)\pi^{2} \equiv 0 \; \text{mod} \; \pi^{3}.$$
For case $p = 2$, this is equivalent to the following
condition\footnote"*"{K. Murasugi informed us of this condition and
examples which are obtained by the relation between the Alexander
polynomial and Milnor invariants [Mu].}
$$
\cases
\text{lk}(K_{1},K_{2}) \equiv 0 \; \text{mod}\; 8, \; \mu(1122) \equiv 0
\; \text{mod} \; 2 \\
\text{or}\\
\text{lk}(K_{1},K_{2}) \equiv 0 \; \text{mod}\; 4,
\text{lk}(K_{1},K_{2}) \not\equiv 0 \; \text{mod}\; 8, \mu(1122) \equiv
1 \; \text{mod}\; 2.
\endcases
\tag 4.2.2.4
$$
For example, the Whitehead link $L = K_{1}\cup K_{2}$ satisfies
$\text{lk}(K_{1},K_{2}) = 0, \mu(1122) = 1$ and so $e_{3} = 1, e_{4} =
0$, hence $H_{1}(M,{\bold Z}_{p}) = \widehat{\Cal O}/{\frak p}^{3}$. For
the 2-bridge link of type (48,37), the latter condition of (4.2.2.4) is
satisfied and so $H_{1}(M,{\bold Z}_{2}) = {\bold Z}/2^{k}{\bold Z}$,
$k\geq 4$.
\bigskip

\head{\bf 5. Iwasawa type formulas for the $p$-homology groups of
$p^{m}$-fold cyclic branched covers} \endhead
\medskip
\flushpar
5.1. {\it Asymptotic formula for the $p$-homology of $p^{m}$-fold
covers.} For $m \geq 1$, let $M_{m}$ be the completion of the
$p^{m}$-fold subcover of $X_{\infty}$ over $X_{L}$ so that $M_{m}$ is a
$p^{m}$-fold cyclic cover of $S^{3}$ branched along $L$. In this last
Section, we are concerned with the asymptotic behavior of the order of
$H_{1}(M_{m},{\bold Z}_{p})$  as $m \rightarrow \infty$ using the
standard argument in Iwasawa theory. As in Section 4, we start again
with the following isomorphisms
$$ \gather \widehat{A_{L}}^{red} \; \simeq \; (H_{1}(X_{\infty},{\bold
Z})\otimes_{\Lambda}\widehat{\Lambda})\oplus \widehat{\Lambda}, \tag
5.1.1 \\
H_{1}(M_{m},{\bold Z}_{p}) \; \simeq \; (H_{1}(X_{\infty},{\bold
Z})\otimes_{\Lambda}\widehat{\Lambda})/(\nu_{p^{m}}(1+X)). \tag 5.1.2
\endgather $$
\flushpar From these, we get immediately an extension of a theorem of M.
Dellomo [D] for a knot.
\medskip
\flushpar
{\bf Proposition 5.1.3} {\it For a link $L$ whose $p$-adic Milnor
invariants are all zero, for example a boundary link,  we have
$H_{1}(M_{m},{\bold Z}_{p}) = {\bold Z}_{p}^{(p^{m}-1)(n-1)}$ for $m\geq
1$. In particular, $H_{1}(M_{m},{\bold Z}_{p}) = 0$ for $m \geq 1$ if
$L$ is a knot. }
\medskip
\flushpar
{\it Proof.} In fact, $\widehat{A_{L}}^{red} = \widehat{\Lambda}^{n}$
for such a link $L$ by Corollary 3.2.5. Hence $H_{1}(X_{\infty},{\bold
Z})\otimes_{\Lambda}\widehat{\Lambda} = \widehat{\Lambda}^{n-1}$ by
(5.1.1) and so $H_{1}(M_{m},{\bold Z}_{p}) = {\bold
Z}_{p}^{(p^{m}-1)(n-1)}$ by (5.1.2).$\qed$
\medskip

In the following, we assume $n \geq 2$. By (5.1.1), the $0$-th
elementary ideal $E_{0}(H_{1}(X_{\infty},{\bold
Z})\otimes_{\Lambda}\widehat{\Lambda})$ over $\widehat{\Lambda}$ is same
as the $1$st completed Alexander ideal $\widehat{E_{1}}(L)$ (3.1.5).
Note that the $1$st $p$-adic Alexander series $\widehat{\Delta_{1}}(L)$
is  given as the greatest common divisors of all $n-1$ minors of the
reduced $p$-adic Traldi matrix and so it is written by the form
$$ \widehat{\Delta_{1}}(L) = X^{n-1}\cdot \widehat{\nabla_{L}}$$
where  we call $\widehat{\nabla_{L}}$ the {\it $p$-adic Hosokawa series}
of $L$. Then by (5.1.2), we have the following formula on the order
$|H_{1}(M_{m},{\bold Z}_{p})|$ which is seen as the $p$-primary part of
the well known formula by S. Kinoshita and F. Hosokawa [KT] (See also
[MM]). Here we interpret $|H_{1}(M_{m},{\bold Z}_{p})| = 0$ to mean
$H_{1}(M_{m},{\bold Z}_{p})$ is infinite.
\medskip
\flushpar
{\bf Proposition 5.1.4.}  $\; |H_{1}(M_{m},{\bold Z}_{p})| =
p^{m(n-1)}\displaystyle{\prod  \Sb \zeta^{p^{m}} = 1\\ \zeta\neq 1
\endSb |\widehat{\nabla_{L}}(\zeta -1)|_{p}^{-1}}.$
\medskip

Now, we assume $H_{1}(M_{m},{\bold Z}_{p})$ is finite for any $m$ and
see the asymptotic behaviour of the order $|H_{1}(M_{m},{\bold Z}_{p})|$
 as $m \rightarrow \infty$. For this, we recall the following standard
facts from Iwasawa theory.
We call a polynomial $g(X) \in {\bold Z}_{p}[X]$ {\it distinguished} if
$g(X) = X^{d} + a_{d-1}X^{d-1} + \cdots + a_{0}$, $a_{i} \equiv 0$ mod
$p$ for $0\leq i \leq d-1$.
\medskip
\flushpar
{\bf Lemma 5.1.5} ($p$-adic Weierstrass preparation theorem [W, Theorem
7.3]). {\it A non-zero element $f(X) \in \widehat{\Lambda}$ is written
uniquely as
$$ f(X) = p^{\mu}g(X)u(X)$$
where $\mu$ is a non-negative integer, $g(X)$ is a distinguished
polynomial and $u(X) \in \widehat{\Lambda}^{\times}$}.
\medskip
\flushpar
{\bf Lemma 5.1.6} ([W,Theorem 7.14]). {\it Let $f(X) \in
\widehat{\Lambda}$ and assume $f(\zeta -1) \neq 0$ for any primitive
$p^{m}$-th root $\zeta$ of 1 for $k \geq 1$. Write $f(X) =
p^{\mu}g(X)u(X)$ according to Lemma 5.5  and define $\lambda$ by the
degree of $g(X)$. Then there is an integer $\nu$ independent of $k$ such
that we have the equality}
$$ \text{ord}_{p}(\prod \Sb \zeta^{p^{m}} = 1\\ \zeta\neq 1 \endSb
f(\zeta -1))
= \lambda m + \mu p^{m} + \nu$$
{\it for sufficiently large $m$.}
\medskip
\flushpar
For the convenience of the reader, we include herewith  a proof of Lemma
5.1.6.
\medskip
\flushpar
{\it Proof of Lemma 5.1.6.} Since $u(x) \in \widehat{\Lambda}^{\times}$,
we have
$$ \text{ord}_{p}(\prod \Sb \zeta^{p^{m}} = 1\\ \zeta\neq 1 \endSb
f(\zeta -1)) = (p^{m}-1)\mu + \text{ord}_{p}(\prod \Sb \zeta^{p^{m}} =
1\\ \zeta\neq 1 \endSb g(\zeta -1)).$$
Write $g(X) = X^{\lambda} + a_{\lambda -1}X^{\lambda -1} +\cdots +
a_{0}$, $a_{i} \equiv 0$ mod $p$ $(0\leq i \leq \lambda -1)$. For a
primitive $p^{l}$-th root $\zeta$ of $1$ $(1\leq l \leq k)$, one has the
equality $(p) = (\zeta - 1)^{\phi(p^{l})}$ of ideals of ${\bold
Z}[\zeta]$ and so $\text{ord}_{p}((\zeta - 1)^{\lambda}) =
\frac{\lambda}{\phi(p^{l})}$ where $\phi(x)$ is the Euler function.
Therefore, if $l$ is large enough, $\text{ord}_{p}((\zeta -
1)^{\lambda}) < \text{ord}_{p}(a_{i}(\zeta - 1)^{i})$ for $0\leq i \leq
\lambda - 1$ and so $\text{ord}_{p}(g(\zeta - 1)) =
\text{ord}_{p}((\zeta - 1)^{\lambda})$. Hence, there is a constant $C$
 independent of $k$ such that for sufficiently large $m$, we have
$$ \text{ord}_{p}(\prod \Sb \zeta^{p^{m}} = 1\\ \zeta\neq 1 \endSb
g(\zeta -1)) = \text{ord}_{p}(\prod \Sb \zeta^{p^{m}} = 1\\ \zeta\neq 1
\endSb (\zeta - 1)^{\lambda}) + C = \text{ord}_{p}(p^{m\lambda}) + C =
\lambda m + C. \qed $$
\medskip
\flushpar
To apply Lemma 5.1.6 to $\widehat{\nabla_{L}}$, write
$$\widehat{\nabla_{L}} = p^{\mu (L;p)}g(L;p)u(L;p) $$
where $\mu (L;p)$ is a nonnegative integer, $g(L;p)$ is a distinguished
polynomial and $u(L,p) \in \widehat{\Lambda}^{\times}$ according to
Lemma 5.1.5 and set $\lambda (L;p) = \text{deg}(g(L;p))$. Then
Proposition 5.1.4 and Lemma 5.1.6 yield the following
\medskip
\flushpar
{\bf Theorem 5.1.7.} {\it Notation and assumption being as above, there
is a constant $\nu (L;p)$ depending only on $L$ and $p$ such that we
have }
$$ \text{ord}_{p}(|H_{1}(M_{m},{\bold Z}_{p})|) = (n-1 + \lambda(L;p))m
+ \mu(L;p)p^{m} + \nu(L;p)$$
{\it for sufficiently large $m$.}
\medskip
\flushpar
We call the invariants $\lambda(L;p), \mu(L;p)$ the {\it Iwasawa
$\lambda$, $\mu$-invariants} of $L$ with respect to $p$ respectively
after the model of the Iwasawa invariants in the theory of ${\bold
Z}_{p}$-extensions [Iw].
\medskip
\flushpar
5.2. {\it  Examples.}
\medskip
\flushpar
1. Let $L$ be the Whitehead link. We then have
$$ \widehat{T_{L}} = \pmatrix X_{1}X_{2}^{2} & -X_{1}^{2}X_{2}\\
X_{1}^{2}X_{2} & -X_{1}X_{2}^{2}
\endpmatrix $$
and $\widehat{\nabla_{L}} = X^{2}$. Hence, we have $\lambda(L;p) = 2,
\mu(L;p) = 0$ and $\text{ord}_{p}(H_{1}(M_{m},{\bold Z}_{p})) = 3m$ for
$m\geq 1$.
\medskip
\flushpar
2. Let $L = K_{1}\cup K_{2}\cup K_{3}$ be the Borromean rings so that we
can take
$ y_{1} = [x_{3},x_{2}], y_{2} = [x_{3},x_{1}], y_{3} = [x_{1},x_{2}]$.
Then we can compute all Milnor number needed to get the reduced Traldi
matrix
$$ \widehat{T_{L}}^{red} = \pmatrix
X + X^{2} & -X^{2} & -X\\
-X & X+X^{2}& -X^{2}\\
-X^{2} & -X & X+X^{2}
\endpmatrix$$
and so $\widehat{\nabla_{L}} = 1+X+X^{2}$. Hence, we have $\lambda(L;p)
= \mu(L;p) = 0$ and $\text{ord}_{p}(H_{1}(M_{m},{\bold Z}_{p})) = 2m$
 for $m\geq 1$.
\medskip
\flushpar

 {\eightrm jonh\@maths.usyd.edu.au; matei\@ms.u-tokyo.ac.jp;
morisita\@kenroku.kanazawa-u.ac.jp}

\enddocument